\documentclass[12pt]{article}
\usepackage{amssymb,amsmath,amsthm,latexsym}
\usepackage[cp1251]{inputenc}
\usepackage{dsfont}
\usepackage{amssymb}
\usepackage{indentfirst}
\usepackage{graphicx}
\usepackage{multicol}
\usepackage{amsfonts, array, hhline}
\usepackage{color}
\begin{document}

\newcommand\blfootnote[1]{%
  \begingroup
  \renewcommand\thefootnote{}\footnote{#1}%
  \addtocounter{footnote}{-1}%
  \endgroup
}

\title{Chromatic numbers of spheres}
\author{Roman Prosanov\thanks{The author is supported by the grant $200021_-169391$ of Swiss National Science Foundation and supported in part by the grant from Russian Federation President Program of Support for Leading Scientific Schools 6760.2018.1.}}
\date{}
\AtEndDocument{\bigskip{\footnotesize%
  \textsc{Universit\'{e} de Fribourg, Chemin du Mus\'{e}e 23, CH-1700 Fribourg, Switzerland} \par
  \textsc{Moscow Institute Of Physics And Technology, Institutskiy per. 9, 141700, Dolgoprudny, Russia} \par
  \textit{E-mail}: \texttt{rprosanov@mail.ru}
}}
\maketitle

\renewcommand{\refname}{Bibliography}
\renewcommand{\proofname}{Proof}
\renewcommand{\figurename}{Figure}
\renewcommand{\le}{\leqslant}
\renewcommand{\leq}{\leqslant}
\renewcommand{\ge}{\geqslant}
\renewcommand{\geq}{\geqslant}
\renewcommand{\mathds}{\mathbb}
\newcommand{\e}{\varepsilon}
\newcommand{\R}{\mathbb{R}}

\begin{abstract}
The chromatic number of a subset of Euclidean space is the minimal number of colors sufficient for coloring all points of this subset in such a way that any two points at the distance 1 have different colors. We give new upper bounds on chromatic numbers of spheres. This also allows us to give new upper bounds on chromatic numbers of any bounded subsets.
\end{abstract}


\section{Introduction}

The \textit{chromatic number} of the Euclidean space $\chi(\R^n)$ is the minimal number of colors sufficient for coloring all points of $\R^n$ in such a way that any two points at the distance 1 have different colors. Determining $\chi(\R^n)$ is considered as an important problem of discrete geometry. For the history and the overview of this problem see~\cite{Ra2}, \cite{Rai2}, \cite{Rai3}, \cite{Ra}, \cite{So}.

Even in the case $n=2$ the exact value of the chromatic number is unknown. The best current bounds are $$5 \leq \chi(\mathds{R}^2) \leq 7,$$ The lower bound is due to de Gray~\cite{dG}, the upper is due to Isbell~\cite{So}.

In the case of arbitrary $n$ Raigorodskii (the lower bound, \cite{Ra1}), Larman and Rogers (the upper bound, \cite{Lar}) proved that $$(1.239 + o(1))^n \leq \chi(\mathds{R}^n) \leq (3 + o(1))^n.$$

For small values of $n$ better lower bounds are known. For the overview of the recent progress we refer to~\cite{CR}.

The chromatic number may be defined for an arbitrary metric space (see for example~\cite{Ku1}, \cite{Ku2}). It also can be defined with any positive real number instead of 1 in the definition of the chromatic number (we call this number the \emph{forbidden distance}). The space $\R^n$ admits homothety, therefore, $\chi(\R^n)$ does not depend on the choice of the forbidden distance, but in general case the chromatic number essentially depends on it.

We consider the case of a spherical space. Let $\chi(S^n_R)$ be the minimal number of colors needed to color the Euclidean sphere $S^n_R$ of radius $R$ in $\R^{n+1}$ in such a way that any two points of $S^n_R$ at the Euclidean distance 1 have different colors. It is clear that $\chi(S^n_{1/2}) = 2$ and $\chi(S^n_R)=1$ for $R < 1/2$. Note that this is the case, when the chromatic number depends on the forbidden distance. But the problem of determining the chromatic number of $S^n_R$ with a forbidden distance $d$ is equivalent to determining the chromatic number of $S^n_{R/d}$ with the forbidden distance 1.


In 1981 Erd\H{o}s conjectured that for any fixed $R > 1/2$, $\chi(S^n_R)$ is growing as $n$ tends to infinity. In 1983 this was proved by Lov\'asz \cite{Lo1} using an interesting mixture of combinatorial and topological techniques. Among other things, in this paper Lov\'asz claimed that for $R<\sqrt{\frac{n+1}{2n+4}}$ (i.e. when the side of regular $(n+1)$-dimensional simplex inscribed in our sphere is less than 1) we have $\chi(S^{n}_R) = n+1$. However, in 2012 Raigorodskii~\cite{Rai1}, \cite{Rai} showed that this statement is wrong. In \cite{Rai} it was shown that actually for any fixed $R > 1/2$ the quantity $\chi(S^n_R)$ is growing exponentially. Some improvements of lower bounds were obtained in \cite{Ko}, \cite{KRa}.

It is clear that $$\chi(S^n_R) \leq (3+o(1))^n$$ because $S^n_R$ is a subset of $\R^{n+1}$. Despite the remarkable interest to this problem there are no better upper bounds in general. For spheres of small radii ($R < 3/2$) the work of Rogers \cite{Ro2} easily implies a much stronger bound. Consider a spherical cap on $S^n_R$ of such radius that the Euclidean diameter of this cap is less than 1. Then we cover $S^n_R$ with copies of this cap and paint every cap in its own color. This establishes the bound $$\chi(S^n_R) \leq (2R + o(1))^n.$$

In this paper we prove a new upper bound on $\chi(S^n_R)$ in the case of $R > \frac{\sqrt 5}{2}$. More precisely, define

$$x(R) = \left\{
\begin{aligned}
& \sqrt{5 - \frac{2}{R^2} + 4\sqrt{1-\frac{5R^2-1}{4R^4}}} , &\qquad R > \frac{\sqrt 5}{2} &  \\
& 2R, &\qquad \frac{1}{2} < R \leq \frac{\sqrt 5}{2} &
\end{aligned}
\right. $$

\vskip+0.2cm	
{\bfseries Theorem 1.} \textit {For $R > \frac{1}{2}$ we have $\chi(S^n_R) \leq (x(R)+o(1))^n.$}
\vskip+0.2cm

It is clear that the base of exponent is always less than 3. (However, it tends to 3 as $R$ tends to infinity.) Further, it will be evident that it is less than $2R$ over the interval $\left(\frac{\sqrt 5}{2}; \frac{3}{2}\right)$. Thus, we improve the current bounds for all $R$ that is not in the interval $\frac{1}{2} < R \leq \frac{\sqrt 5}{2}$. In the latter case the method of our proof breaks down, but it provides another proof of the bound $\chi(S^n_R) \leq (2R + o(1))^n$.

Let $B^{n+1}_R \subset \R^{n+1}$ be a Euclidean ball of radius $R$ (centered in the origin). By $\chi(B^{n+1}_R)$ denote the chromatic number of $B^{n+1}$ (with forbidden distance 1). The construction in the proof of Theorem 1 also implies the following

\vskip+0.2cm	
{\bfseries Theorem 2.} \textit {For $R > \frac{1}{2}$ we have $\chi(B^{n+1}_R) \leq (x(R)+o(1))^n.$}
\vskip+0.2cm

The Erd\H{o}s--de Bruijn theorem \cite{EdB} states that the chromatic number of the Euclidean space $\R^n$ is reached at some finite distance graph embedded in this space. Hence, Theorem 2 connects $\chi(\R^n)$ with the radius of circumscribed sphere of this graph. The author is grateful to A.B. Kupavskii for the remark that Theorem 1 should imply Theorem 2.

It is of interest to mention the problem of determining \emph{the measurable chromatic number} $\chi_m(\R^n)$, which is defined in the same way, but with the extra condition that all monocolored sets are required to be measurable. In this case upper bounds remain to be the same, but additional analytic techniques can be applied to establish better lower bounds. Thus, in~\cite{BPT} in was proved that $$(1.268-o(1))^n \leq \chi_m(\R^n).$$ The best lower bound for $n=2$ was obtained by Falconer~\cite{Fa}, the best lower bounds for some other small values of $n$ can be found in~\cite{DfW}.

Another fruitful area of research is to consider colorings with more restrictions on a monocolored set. Some results in this direction were obtained in~\cite{Be}, \cite{Be2}, \cite{BeRa}, \cite{PoRa}, \cite{Pr2}, \cite{Sa2}, \cite{Sa}.

This paper is organized as follows. The next section contains a brief exposition of our technique and provides further references to the bibliography. In Section 3 we set up notation and terminology. Section 4 is devoted to the proof of Theorem 1 and Section 5 provides the implication of Theorem 2.

\section{Summary of the technique}

Generally, in our proof of Theorem 1 we follow the approach of Larman and Rogers in~\cite{Lar} about the chromatic number of the Euclidean space. We construct some set on $S^n_R$ without a pair of points at the distance 1, bound its density and then cover the whole sphere with copies of this set. But the realization of every item of this plan in \cite{Lar} can not be generalized directly to the spherical case. For instance, the construction of the set without distance 1 is strongly based on theory of lattices in $\R^n$. Therefore, we should provide some new ideas for our case.

During our proof we should turn to geometric covering problems. We need to cover a sphere with copies of some disconnected set. There is a well-developed theory of economical coverings. The most common tool is the Rogers theorem \cite{Ro1} (and its relatives) on periodical coverings of a Euclidean space with copies of a convex body. There are different approaches to prove results of that type. Most of proofs essentially use the convexity of the body and special properties of Euclidean spaces. This is not appropriate in our situation.

A new approach to geometric covering problems was proposed by Nasz\'odi in~\cite{Na}. He suggested to construct some finite hypergraph in such a way that an edge-covering of this hypergraph implies a desired geometric covering. The covering numbers of finite hypergraphs can be studied via the famous result on the ratio between the optimal covering number and the optimal fractional covering number (see the next section for the definition of a fractional covering and the statement of this lemma). Surprisingly, the fractional covering number of our hypergraph can be bounded from geometric setting. In \cite{Na} these ideas were used to provide new proofs of the Rogers theorem on periodical coverings of $\R^n$ and the Rogers result on coverings of a sphere with caps. In~\cite{NP} Naszd\'odi and Polyanskii proved new results on geometric multi-coverings with the help of this approach.

The work \cite{Na} pursues its own purposes and its results can not be straightly applied in our setting. Therefore, in section 4.3 we present in details a self-contained proof of the covering part and adapt Nasz\'odi's arguments for our needs. The author already used similar techniques in the papers \cite{Pr1} and \cite{Pr2} for analogous problems in the Euclidean space. 

It should be mentioned that the concept of  ``fractional'' geometric coverings was also suggested by Artstein-Avidan, Raz and Slomka in the papers \cite{AR} and \cite{AS}. For a deeper discussion of geometric coverings we refer the reader to the recent survey \cite{NaS} by Nasz\'odi.

\section{Preliminaries}

\begin{figure}
\begin{center}
\includegraphics[scale=0.2]{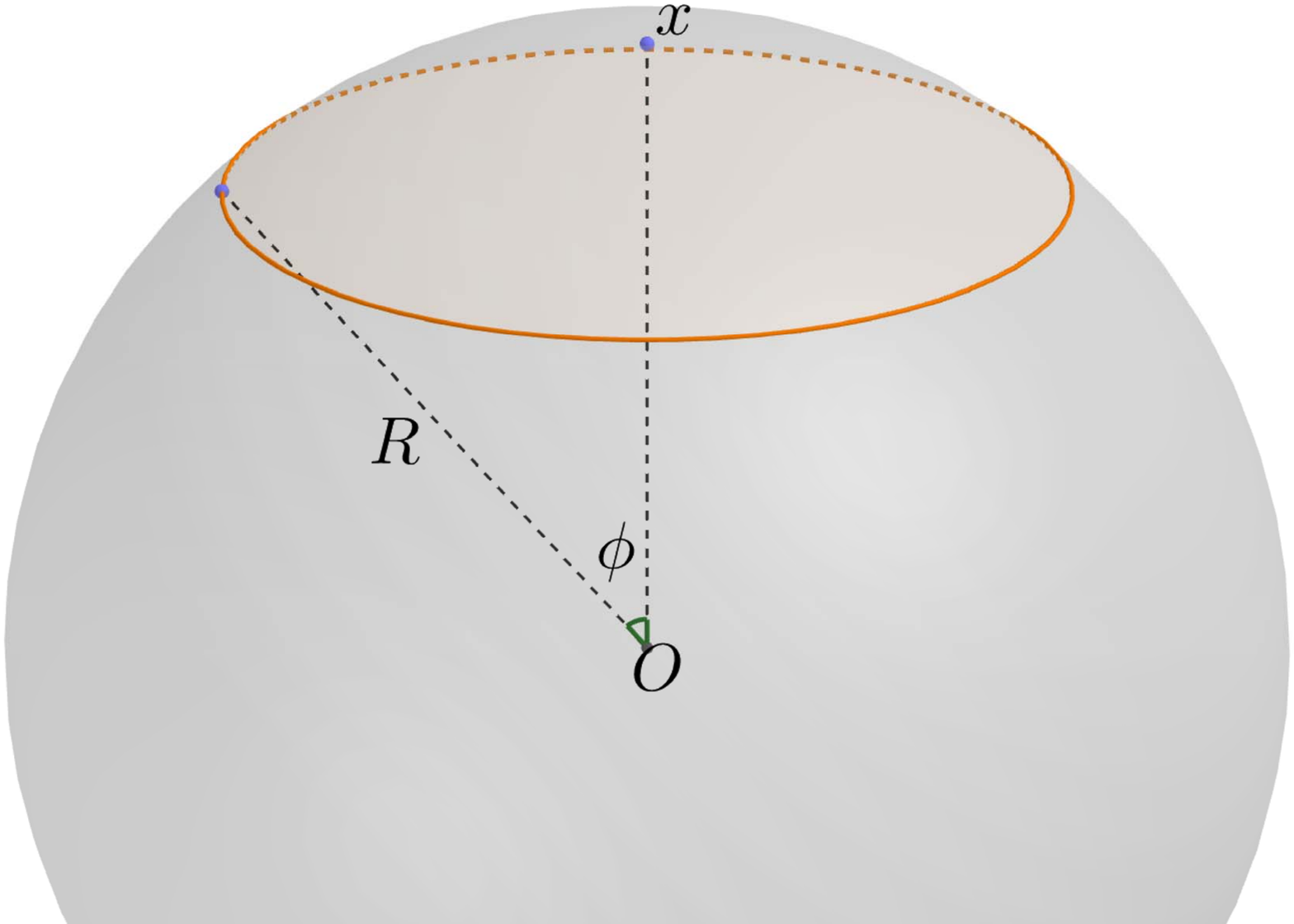}
\caption{A spherical cap of angular radius $\phi$.}
\label{cap}
\end{center}
\end{figure}

Let $C(x, \phi)$ be a spherical cap (Figure~\ref{cap}) with center $x \in S^n_R$ of angular radius $\phi \leq \pi/2$, i.e. the set of all points $x' \in S^n_R$ such that the angle $xox' \leq \phi$. The Euclidean diameter of $C(x, \phi)$ is $2R\sin\phi$.

Denote be ${\rm vol}(Z)$ the spherical volume of a measurable set $Z \subset S^n_R$ and by $\rho(.)$ the usual probability measure on $S^n_R$, i.e. $\rho(Z)=\frac{{\rm vol}(Z)}{{\rm vol}(S^n_R)}$. We will use the word ``density'' for $\rho(.)$.

Let $B(o, r)$ be the Euclidean ball in $\R^{n+1}$ with center $o$ and radius $r$.

Define $\Theta(\phi) = \rho(C(x, \phi))$. The following auxiliary result will be needed in Section 4.3.

\vskip+0.2cm	
{\bfseries Lemma 1, B\"{o}r\"{o}czky--Wintsche, \cite{BW}.} \textit {Let ${0<\phi<\pi/2}$ and $1<t<\frac{\pi}{2\phi}$. Then
$$\Theta(t\phi) < t^n\Theta(\phi),$$
$$\Theta(\phi) > \frac{\sin^n{\phi}}{\sqrt{2\pi(n+1)}}.$$}
\vskip+0.2cm

Consider an arbitrary set $X$ and a family of subsets $\mathcal{F} \subset 2^X$. Denote by $\tau(X, \mathcal{F})$ the minimal cardinality of a subfamily $\mathcal{G} \subseteq \mathcal{F}$ such that the union of all subsets in $\mathcal{G}$ contains $X$ (if the union of all subsets in $\mathcal{F}$ does not contain $X$, then set $\tau(X, \mathcal{F})$ to be equal $\infty$).

The following special case is particulary important. Let $G = (\mathcal{V}, \mathcal{E})$ be a finite hypergraph. A \textit{covering} of this hypergraph is a family of edges such that their union contains all vertices. \textit{The covering number} $\tau(G)$ is the minimal size of a covering of $G$. In the notation of the previous paragraph, $\tau(G) = \tau(\mathcal V, \mathcal E)$.

A \textit{fractional covering} of $G$ is a function $\nu: \mathcal{E} \rightarrow [0; +\infty)$ such that for every $v \in \mathcal{V}$, we have $\sum\limits_{E \in \mathcal{E}: v \in E} \nu(E) \geq 1$. Define the \textit{fractional covering number}

\vskip+0.2cm
$\tau^*(G)= \inf \{ \sum\limits_{E \in \mathcal{E}}\nu(E): \nu$ is a fractional covering of $G \}$.
\vskip+0.2cm

The following lemma establishes a connection between the integral covering number and the fractional one.

\vskip+0.2cm
{\bfseries Lemma 2, \cite{Jo, Lo, St}.} \textit {Let $G$ be a finite hypergraph. Then
$$
\tau(G) < \left(1+ \ln\left(\underset{E \in \mathcal{E}}{\max} (|E|)\right)\right)\tau^*(G).$$
}
\vskip+0.2cm

\section{Proof of Theorem 1}

\subsection{Construction of a set without distance 1}

Let $\phi < \pi /4$ be a fixed angle. Consider a set $X \subset S^n_R$ of maximal cardinality such that for all $x_1, x_2 \in X$ we have $C(x_1, \phi) \cap C(x_2, \phi) = \emptyset$. Then $\bigcup\limits_{x \in X} C(x, 2\phi) = S^n_R$. Indeed, if there is $y \notin \bigcup\limits_{x \in X} C(x, 2\phi)$, then $C(y, \phi)$ does not intersect any cap $C(x, \phi)$ for $x \in X$. This contradicts the maximal cardinality of $X$.

Consider the Voronoi tiling $\Psi$ of the sphere $S^n_R$ corresponding to the set $X$. By $\psi_x$ denote the unique spherical convex polytope of this tiling containing a point $x \in X$. Note that $C(x, \phi) \subset \psi_x \subset C(x, 2\phi)$.

Let $L$ be the tangent hyperplane to $S^n_R$ at a point $x \in X$. Define a map $$f_{x,\lambda}: C(x, \pi/2) \rightarrow C(x, \pi/2).$$ Take an arbitrary point $a \in C(x, \pi/2)$. Let $p$ be the line through $a$ orthogonal to $L$. Consider the homothety of $p$ with center $x$ and coefficient $0 < \lambda < 1$. Its image intersects the half-sphere $C(x, \pi/2)$ at a point $a'$. We define $f_{x, \lambda} (a)$ to be $a'$. Write $\psi'_x$ for $f_{x, \lambda} (\psi_x)$ and $\Psi'$ for $\bigcup\limits_{x \in X} \psi'_x$.

Determine $\gamma$ by the equation $\sin\gamma=\lambda\sin(2\phi)$. It is clear that $\psi'_x$ is contained in the spherical cap $C(x, \gamma)$.  Therefore, the Euclidean diameter of $\psi'_x$ is not greater than $2R\lambda\sin(2\phi)$.

We want to give a lower bound on the minimal distance between two sets $\psi'_x$ and $\psi'_y$. To this purpose we need the following technical proposition. We postpone its proof until Subsection 4.5.

\vskip+0.2cm
{\bfseries Proposition 1.} \textit{Let the angle $\alpha$ be determined by the equation $\sin \alpha = \lambda \sin \phi$. Then the angular distance between $\psi'_x$ and the boundary of $\psi_x$ is not less than $\phi-\alpha$. }
\vskip+0.2cm

Proposition 1 implies that the angular distance between the sets $\psi'_x$ and $\psi'_y$ ($x,y \in X$) is not less than $2(\phi-\alpha)$. Indeed, consider a pair of points $z_x \in \psi'_x$ and $z_y \in \psi'_y$ and the unique shortest geodesic arc between these points. There are two disjoint sub-arcs: from $z_x$ to the boundary of $\psi_x$ and from the boundary of $\psi_y$ to $z_y$. Each of these sub-arcs has the angular length not less than $\phi-\alpha$. Therefore, the angular distance between $z_x$ and $z_y$ is not less than $2(\phi-\alpha)$ and the Euclidean distance between the sets $\psi'_x$ and $\psi'_y$ is not less than $2R\sin(\phi-\alpha)$.

Next, we determine $\phi$ and $\lambda_0$ such that

\begin{equation}\label{se}
\left\{
\begin{aligned}
2R\lambda_0\sin(2\phi)=1\\
\sin{\alpha}=\lambda_0\sin{\phi} \\
2R\sin(\phi-\alpha)=1\\
\end{aligned}
\right.
\end{equation}

All previous observations show that for $\lambda < \lambda_0$ the set $\Psi'$ does not contain a pair of points at the distance 1.

To solve this equation system we first determine $\lambda_0$ as a function of $\phi$. Comparing the left sides of the first and the third equations we obtain
$$2\lambda_0\sin \phi \cos \phi = \sin \phi \cos \alpha - \sin \alpha \cos \phi=$$
$$=\sin \phi \sqrt{1-\lambda_0^2\sin^2 \phi} - \lambda_0\sin \phi \cos \phi.$$
Therefore,
$$2\lambda_0 \cos \phi = \sqrt{1-\lambda_0^2\sin^2 \phi} -\lambda_0\cos \phi,$$
$$9\lambda_0^2\cos^2 \phi = 1 - \lambda_0^2 \sin^2 \phi,$$
$$\lambda_0^2 = \frac{1}{1+8\cos^2 \phi}.$$

We substitute this expression in the first equation and get the equation for $\phi$:
$$1+8\cos^2 \phi = 16R^2\sin^2 \phi \cos^2\phi.$$

Solving this equation we obtain:
$$\cos^2 \phi = \frac{1}{2} - \frac{1}{4R^2} \pm \sqrt{\frac{1}{4}-\frac{5R^2-1}{16R^4}}.$$

Recall that we have the restriction $\phi < \pi/4$. Therefore, we should use only the root with plus sign and should determine, when it is greater than 1/2. We get the inequality
$$\frac{1}{2} - \frac{1}{4R^2} + \sqrt{\frac{1}{4}-\frac{5R^2-1}{16R^4}} > \frac{1}{2}.$$

Solving it we obtain the restriction $R > \frac{\sqrt{5}}{2}$. During the next two subsections we will consider only this case. The latter case $\frac{1}{2} < R \leq \frac{\sqrt{5}}{2}$ will be discussed in Subsection 4.4.

It is obvious that $\lambda_0 > \frac{1}{3}$. Moreover, from the first equation we have $$\lambda_0=\frac{1}{2R\sin{2\phi}}.$$ Therefore, $\lambda_0 > \frac{1}{2R}$.

\subsection{Density estimating}

Consider $x \in X$ and choose an orthogonal coordinate system in $\R^{n+1}$ in such a way that the axis $Ox_{n+1}$ coincides with the ray $Ox$. The sets $\tilde{\Omega}$ and $\tilde{\Omega}'$ are the orthogonal projections of $\psi_x$ and $\psi_{x}'$ to the hyperplane $Ox_1...x_n$. It is easily seen that $\frac{1}{\lambda}\tilde{\Omega}' = \tilde{\Omega}$.

We have

\begin{equation}\label{int}
{\rm vol}(\psi_x) = \underset{\tilde{\Omega}}{\int \ldots \int} \frac{R dx_1 \ldots dx_n}{\sqrt{R^2 - x_1^2 - \ldots - x_n^2}} = \underset{\frac{1}{\lambda}\tilde{\Omega}'}{\int \ldots \int} \frac{R dx_1 \ldots dx_n}{\sqrt{R^2 - x_1^2 - \ldots - x_n^2}} =$$ $$= \underset{\tilde{\Omega}'}{\int \ldots \int} \frac{\lambda^{-n}R dx_1 \ldots dx_n}{\sqrt{R^2 - \frac{x_1^2 + \ldots + x_n^2}{\lambda^2}}}.
\end{equation}

Remind that $\psi_x'$ is contained in $C(x, \gamma)$. Hence, $\tilde{\Omega}' \subset B(O, R\sin\gamma)=B(O, R\lambda\sin(2\phi))$.

Consider the function $$f(r) = \sqrt{\frac{R^2-r^2}{R^2-\frac{r^2}{\lambda^2}}}.$$ This function increases monotonically over $[0; \lambda R)$. Therefore, $f(r)$ reaches its maximal value over the segment $[0; R\lambda\sin(2\phi)]$ at the right endpoint. This enables us to get for an arbitrary $r \in [0; R\lambda\sin(2\phi)]$ the inequality
$$\frac{1}{\sqrt{R^2-\frac{r^2}{\lambda^2}}} \leq \sqrt{\frac{R^2 - R^2\lambda^2\sin^2(2\phi)}{R^2-\frac{R^2\lambda^2\sin^2(2\phi)}{\lambda^2}}}\cdot \frac{1}{\sqrt{R^2-r^2}}.$$

We can now bound the value of the integral in (\ref{int})
$${\rm vol}(\psi_x) \leq \lambda^{-n} \sqrt{\frac{R^2 - R^2\lambda^2\sin^2(2\phi)}{R^2-\frac{R^2\lambda^2\sin^2(2\phi)}{\lambda^2}}} \underset{\tilde{\Omega}'}{\int \ldots \int} \frac{R dx_1 \ldots dx_n}{\sqrt{R^2 - x_1^2 + \ldots + x_n^2}} =$$ $$= \lambda^{-n} \sqrt{\frac{1 - \lambda^2\sin^2(2\phi)}{1-\sin^2(2\phi)}} \cdot {\rm vol}(\psi_x').$$

Hence, we have the following bound on the density

\begin{equation}\label{den}
\rho(\Psi') \geq \min\limits_{x \in X} \frac{{\rm vol}(\psi_x')}{{\rm vol}(\psi_x)} \geq \lambda^{n} \sqrt{\frac{1-\sin^2(2\phi)}{1 - \lambda^2\sin^2(2\phi)}}.
\end{equation}

\subsection{Covering of $S^n_R$ with copies of the set $\Psi'$}

Let $SO(n+1)$ be the group of orientation-preserving isometries of $S^n_R$ (naturally identified with the group of orthogonal $(n+1)\times(n+1)$ matrices of determinant 1). Fix $0< \delta < 1$. By $\psi''_x$ denote $f_{x, (1-\delta)\lambda}(\psi_x)$ and by $\Psi''$ denote the union $\bigcup\limits_{x \in X} \psi''_x$. Let $\mathcal{F'}$ and $\mathcal{F''}$ be the families of all possible images of the sets $\Psi'$ and $\Psi''$, respectively, under the action of $SO(n+1)$.

Determine $\beta$ by the equation $\sin{2\beta}=\lambda\delta\sin{\phi}$. Consider a set $W \subset S^n_R$ of maximal cardinality such that for every $w_1, w_2 \in W$ we have $C(w_1, \beta) \cap C(w_2, \beta) = \emptyset$. It is clear that $S^n_R = \bigcup\limits_{w \in W} C(w, 2\beta)$ (for details see the beginning of Subsection 4.1). 

\vskip+0.2cm
{\bfseries Proposition 2.} $\tau(S^n_R,\mathcal{F'}) \leq \tau(W,\mathcal{F''})$. (The definition of $\tau$ is in Section~3.)
\vskip+0.2cm

\begin{proof}

Consider $\mathcal{A} \subset SO(n+1)$ such that $W \subset \bigcup\limits_{A \in \mathcal{A}} A\Psi''$. It is sufficient to show that $S^n_R \subseteq \bigcup\limits_{A \in \mathcal{A}} A\Psi'$.

Take $w \in \Psi''$. Recall that $S^n_R = \bigcup\limits_{w \in W} C(w, 2\beta)$. If we prove that $C(w, 2\beta) \subset \Psi'$, then it will be evident that $S^n_R = \bigcup\limits_{w \in W} C(w, 2\beta) = \bigcup\limits_{A \in \mathcal{A}} A\Psi'$.

Choose $x \in X$ such that $w \in \psi''_x$. We remind that $L$ is the tangent plane to $S^n_R$ at the point $x$. Define a cylinder $P_x$ as the union of all lines orthogonal to $L$ passed through points of $\psi_x$. Note that it is convex. Indeed, consider a spherical triangle $abx$ ($a, b \in \psi_x \subset C(x,\pi/2)$). Under these assumptions the orthogonal projection of this triangle to $L$ is a convex set (it is enough to examine the planar case only). Therefore, $ab \subset P_x$. By $P'_x$ denote the image of $P_x$ under the homothety of this cylinder with center $x$ and coefficient $\lambda$. Analogously, by $P''_x$ denote the image of $P_x$ under the homothety with center $x$ and coefficient $(1-\delta)\lambda$. We see that $w \in P''_x$. At last, consider another homothety with center $x$ and coefficient $\lambda\delta$. By $Q$ denote the image of $P_x$ under the composition of this homothety and the translation by the vector $w$.

The cylinder $P_x$ contains the cap $C(x, \phi)$. This implies that $P_x$ contains the Euclidean ball $B(x, R\sin{\phi})$. Therefore, $Q$ contains the Euclidean ball $B(w, \lambda\delta R\sin{\phi}) = B(w, R\sin{2\beta})$ and $C(w, 2\beta) \subset (Q \cap S^n_R)$. The proof is completed by showing that $P'_x \supset Q$.

Indeed, $Q = \delta P'_x + w \subset \delta P'_x + P''_x = \delta P'_x + (1-\delta)P'_x = P'_x$.

\end{proof}

Let $G$ be a hypergraph on vertex set $W$ and edge set $\mathcal E $ defined as follows: $$\mathcal{E} = \{W \cap A\Psi' : A \in SO(n+1) \}.$$ 

It follows from the definition that $\tau(W,\mathcal{F''}) = \tau(G)$.

Lemma 2 together with Proposition 2 imply that $$\tau(S^n_R,\mathcal{F'}) < \left(1+ \ln\left(\underset{E \in \mathcal{E}}{\max} (|E|)\right)\right)\tau^*(G).$$

Consider the Haar measure $\sigma$ on $SO(n+1)$ with the normalization condition $\sigma(SO(n+1)) = 1$. For every $E \in \mathcal{E}$ the set $$S(E) = \{A \in SO(n+1): A\Psi'' \cap W = E \}$$ is closed and, therefore, measurable. Define a function $\nu: \mathcal{E} \rightarrow [0; +\infty)$: $$\nu(E) = \frac{\sigma(S(E))}{\rho(\Psi'')}.$$

For every $w \in W$ the set $$S(w) = \{A \in SO(n+1): w \in A\Psi'' \}$$ is measurable too. It is clear that $$\sum\limits_{E \in \mathcal{E}: w \in E} \nu(E) = \frac{\sigma(S(w))}{\rho(\Psi'')}.$$ Indeed, for every $A$ there is a unique edge $E \in \mathcal{E}$ such that $A \in S(E)$. If $A \in S(w)$, then $w \in E$ and $S(E) \subset S(w)$. Moreover, for every  $w \in S^n_R$ we have $$\sigma(S(w))=\sigma(\{A \in SO(n+1): w \in A\Psi'' \}) =$$ $$= \sigma(\{A \in SO(n+1): w \in A^{-1}\Psi'' \} =$$ $$= \sigma(\{A \in SO(n+1): A w \in \Psi'' \} = \rho(\Psi'').$$

For every $w \in W$ we obtain that $\frac{\sigma(S(w))}{\rho(\Psi'')}=1$. Then $\nu$ is a fractional covering of a hypergraph $G$. Therefore, we have $$\tau^*(G) \leq \sum\limits_{E \in \mathcal{E}}\nu(E) = \frac{\sigma(SO(n+1))}{\rho(\Psi'')} = \frac{1}{\rho(\Psi'')}.$$

Substituting $(1-\delta)\lambda$ instead of $\lambda$ into (\ref{den}) we can bound the density $\rho(\Psi'')$:

$$\rho(\Psi'') \geq \lambda^{n}(1-\delta)^n \sqrt{\frac{1-\sin^2(2\phi)}{1 - \lambda^2(1-\delta)^2\sin^2(2\phi)}}.$$

We next give the bound on $$\underset{E \in \mathcal{E}}{\max} (|E|) = \underset{A \in SO(n+1)}{\max} (|W \cap A\Psi''|).$$

It is clear that $|X| \leq \frac{1}{\Theta(\phi)}$ (the definition of $\Theta$ is in Section 3). Determine $\gamma'$ by the equation $\sin\gamma'=\lambda(1-\delta)\sin(2\phi)$. Since $\psi''_x \subset C(x, \gamma')$, it follows that $w \in \psi''_x$ implies $C(w, \beta) \subset C(x, \gamma'+\beta)$. For all $w_1, w_2 \in W$ it is true that $C(w_1, \beta) \cap C(w_2, \beta) = \emptyset$. We can proceed with the bounds on the volumes. We apply Lemma 1 and obtain
$$|W \cap A\Psi''| \leq \frac{|X|\Theta(\gamma'+\beta)}{\Theta(\beta)} < \frac{\left(1 + \frac{\gamma'}{\beta}\right)^n}{\Theta(\phi)} < \left(1 + \frac{\gamma'}{\beta}\right)^n \frac{\sqrt{2\pi(n+1)}}{\sin^n{\phi}}.$$

Finally, we have:
\begin{equation}\label{pre}
\tau(S^n_R,\mathcal{F'}) \leq
\end{equation}

$$\leq \lambda^{-n}(1-\delta)^{-n} \sqrt{\frac{1 - \lambda^2(1-\delta)^2\sin^2(2\phi)}{1-\sin^2(2\phi)}} \cdot$$ $$\cdot\left(1 + n\ln(1+\frac{\gamma'}{\beta}) - n\ln(\sin{\phi}) + \frac{1}{2}\ln(2\pi(n+1)) \right).$$

We substitute $\delta = \frac{1}{2n \ln n}$ and use (for large $ n $)
$$
\Big(1 - \frac{1}{2n \ln n}\Big)^{-n} \leq \exp \Big(\frac{1}{\ln n}\Big) \leq 1+\frac{2}{\ln n}.
$$

The expression $\sqrt{\frac{1 - \lambda^2(1-\delta)^2\sin^2(2\phi)}{1-\sin^2(2\phi)}}$ tends to a constant as $n$ tends to infinity when $R$ is fixed. The angle $\gamma'$ tends to a constant too. Only $\beta$ is of order~$\frac{\lambda\sin\phi}{4n\ln n}$. Therefore, in the right side of (\ref{pre}) $\lambda^{-n}$ is multiplied by the factor of order $n\ln n$. Hence, we get $$\tau(S^n_R,\mathcal{F'}) \leq (\lambda^{-1} + o(1))^n.$$

Together with bounds on $\lambda$ from Subsection~4.1 it finishes the proof of Theorem~1.

\subsection{Notes on the case $\frac{1}{2} < R \leq \frac{\sqrt{5}}{2}$}

Consider $\frac{1}{2} < R \leq \frac{\sqrt{5}}{2}$. We can study the system of equations (\ref{se}) and see that for every $0 < \phi < \pi/4$, $\e > 0$ and $\lambda = \frac{1}{2R\sin(2\phi)+\e}$ the set $\Psi'$ does not contain a pair of points at the distance 1. When we choose $\phi$, we should take into account the lower bound (\ref{den}) on the density of $\Psi'$. As $\phi$ tends to $\pi/4$, $\lambda$ increases and tends to $\frac{1}{2R+\e}$. But we can see that the sub-exponential factor in (\ref{den}) tends to zero in this case. However, for any fixed $\phi$ this quickly becomes insignificant as $n$ tends to infinity. We have $$\chi(S^n_R) \leq (2R\sin(2\phi) + \e + o(1))^n.$$ We can choose sequences $\phi_n \rightarrow \pi/4$ and $\e_n \rightarrow 0$. Then, we get another proof of the inequality $$\chi(S^n_R) \leq (2R + o(1))^n.$$

The monocolored set from the Rogers bound is a cap $C(x, \gamma'')$, where $\gamma''$ is determined by the equation $\sin\gamma''=\frac{1}{2R+\e}$ (see Sections~1 and~3). Our monocolored set $\Phi'$ consists of pieces $\psi_x$. Each $\psi_x$ is a subset of $C(x, \gamma) \subset C(x, \gamma'')$. As $\phi_n$ tends to $\pi/4$, the number of pieces (which is equal to the size of the set $X$) depends sub-exponentially on~$n$ (for example, see \cite{Bo}, bounds on packing densities of caps of radius $\pi/4$). In summary, our monochromatic set is the union of pieces, each piece has a density not greater than the density of Rogers monochromatic set and the number of pieces depends sub-exponentially on $n$. Therefore, it is natural that our construction does not permit us to improve the base of exponent in the Rogers bound in this case. It seems that the careful analysis of the size of $X$ and a subtle choice of $\phi_n$ and $\e_n$ may allow us to show that we can obtain a better sub-exponential factor from our construction rather than from the construction with one cap, but this is not the aim of the present paper.  

\subsection{Proof of Proposition 1}

By definition $\psi_x$ is a spherical polytope. Let $l$ be a great $(n-1)$-dimensional hypersphere containing a facet of this polytope. It is sufficient to prove that the angular distance between $\psi'_x$ and $l$ is not less than $\phi-\alpha$.

Let $\phi_1 \geq \phi$ be the angular distance between $x$ and $l$, let $\alpha_1$ be the angle determined by the equation $\sin \alpha_1 = \lambda \sin \phi_1$. We show that $\phi_1-\alpha_1$ is not less than $\phi-\alpha$. To prove this claim, compute the derivative of $\varphi-\arcsin(\lambda\sin\varphi)$ as a function of $\varphi$. It is equal to $$1-\frac{\lambda\cos\varphi}{\sqrt{1-\lambda^2\sin^2\varphi}} > 0.$$

Hence, we only need to show that the angular distance between $\psi'_x$ and $l$ is not less than $\phi_1-\alpha_1$. Consider the locus of points of $S^n_R$ such that the angular distance between any of these points and $l$ is not less than $\phi_1-\alpha_1$. It is clear that this locus is a pair of closed spherical caps. One of them contains $x$. Let $\tilde{l}$ be the boundary sphere of this cap. Then $\tilde{l}$ is parallel to $l$ and the angular distance between them is equal to $\phi_1-\alpha_1$. We show that $\psi'_x$ is contained in the cap bounded by $\tilde{l}$. It is sufficient to prove it only for boundary points of $\psi'_x$. The further proof is divided in two steps. First, we prove it for points in $\partial\psi'_x$ that have pre-images (under $f_{x, \lambda}$) in $l$. Second, we show it for all points in $\partial \psi'_x$.

\begin{figure}
\begin{center}
\includegraphics[scale=0.55]{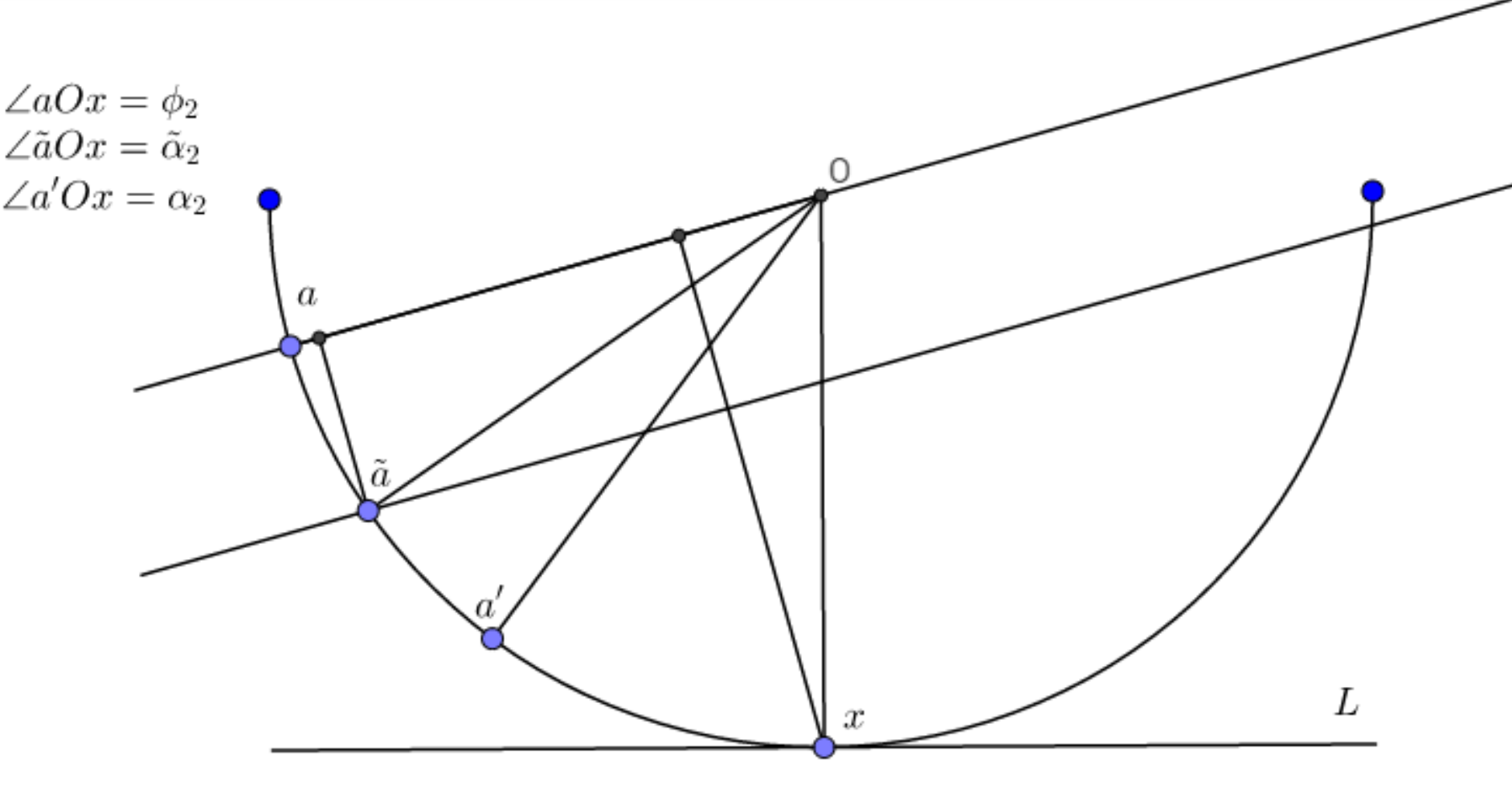}
\caption{The positions of points in Step 1. The orthogonal segments may be helpful to understand the equation (\ref{par}).}
\label{fig1}
\end{center}
\end{figure}

\textit{Step 1.} Consider an arbitrary point $a \in l$. By $a'$ denote $f_{x,\lambda}(a)$ and by $\tilde{a}$ denote the intersection point of $\tilde{l}$ and the arc $xa$. Our goal is to show that the point $\tilde{a}$ lies between the points $a$ and $a'$ (see Figure~\ref{fig1}). Define the angles $\alpha_2 = a'Ox$, $\tilde{\alpha}_2 = \tilde{a}Ox$ and $\phi_2 = aOx \geq \phi_1$. From definitions we see that $$\frac{\sin \alpha_2}{\sin\phi_2} = \lambda = \frac{\sin \alpha_1}{\sin\phi_1}.$$

Consider the point $a_0 \in l$ such that the arc $xa_0$ is orthogonal to $l$. The length of this arc is equal to $\phi_1$ and the length of the sub-arc between $l$ and $\tilde{l}$ is equal to $\phi_1 - \alpha_1$. Look at the two perpendiculars $xy_0$ and $xy_1$, where $y_0$ lies at the segment $Oa_0$ and $y_1$ lies at the segment $Oa$. The hyperspheres $l$ and $\tilde{l}$ are parallel, i.e. the planes in $\R^{n+1}$ determining them are parallel. Therefore, the second plane divides $xy_0$ and $xy_1$ in the same ratio. Using it we conclude the equation

\begin{equation}\label{par}
\frac{\sin (\phi_2-\tilde{\alpha}_2)}{\sin\phi_2} = \frac{\sin(\phi_1-\alpha_1)}{\sin{\phi_1}}.
\end{equation}

It suffices to show that $$\sin \tilde{\alpha}_2 \geq \sin \alpha_2 = \lambda \sin \phi_2 = \frac{\sin \alpha_1 \sin \phi_2}{\sin \phi_1}.$$

It is equivalent to $$\frac{\sin \tilde{\alpha}_2}{\sin\phi_2} \geq \frac{\sin \alpha_1}{\sin\phi_1}.$$

Define $$c = \frac{\sin (\phi_2-\tilde{\alpha}_2)}{\sin\phi_2} = \frac{\sin(\phi_1-\alpha_1)}{\sin{\phi_1}} < 1.$$

Consider $\alpha(\varphi) = \varphi - \arcsin (c\cdot\sin \varphi)$ and $$f(\varphi) = \frac{\sin \alpha (\varphi)}{\sin\varphi} = \cos(\arcsin(c\cdot\sin\varphi)) - c\cdot\cos\varphi.$$
$$f'(\varphi)=c\cdot\sin\varphi-c^2\cdot\sin\varphi\cos\varphi\frac{1}{\sqrt{1-c^2\sin^2\varphi}}.$$

For $0<\varphi < \pi/2$ the inequality $f'(\varphi)>0$ is equivalent to $$1 > c^2(\sin^2 \varphi + \cos^2 \varphi) = c^2.$$

Therefore, we obtain that this function is growing over the interval $(0; \pi/2)$. This finishes the proof of Step 1.

\begin{figure}
\begin{center}
\includegraphics[scale=0.18]{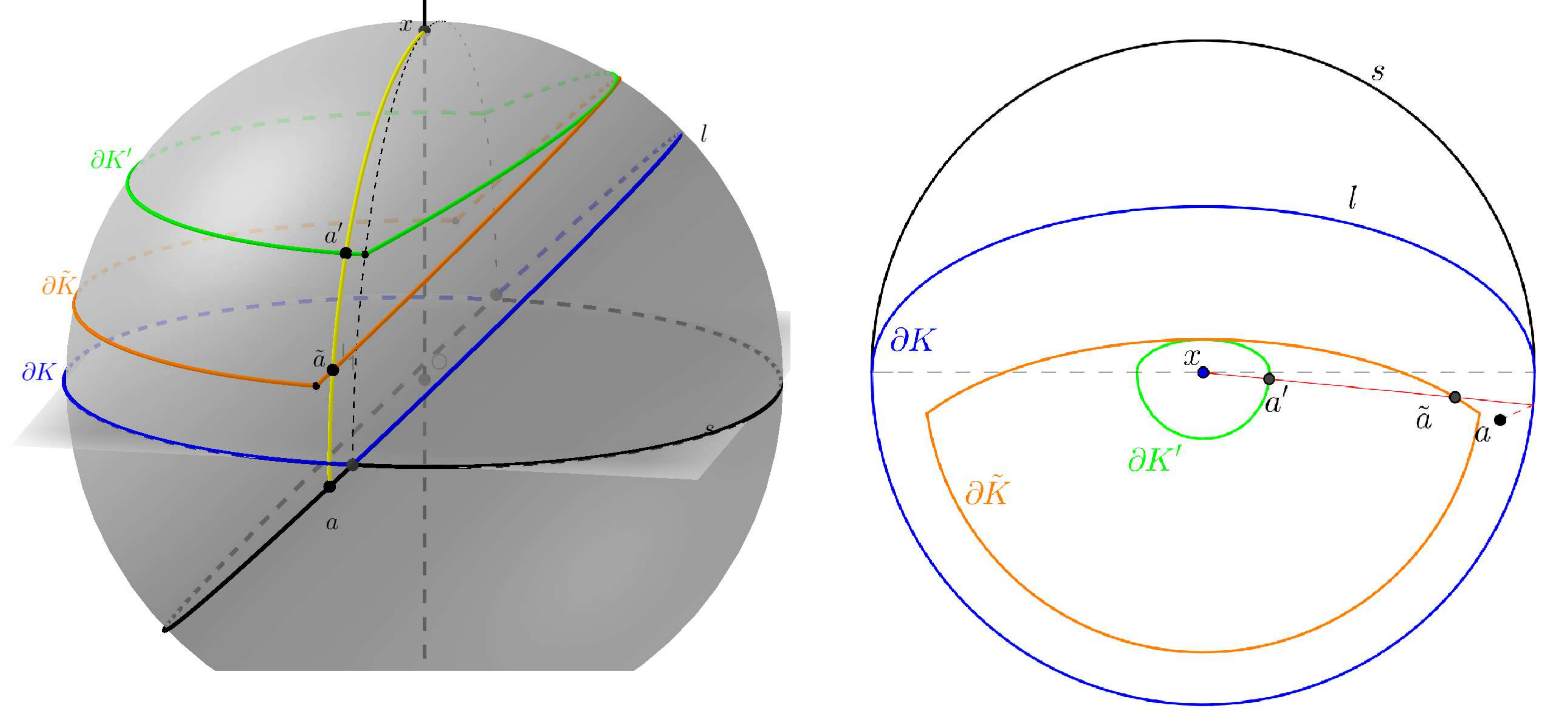}
\caption{The positions of the sets $K$, $\tilde{K}$ and $K'$ in the case $n=2$. Left: on the sphere. The arcs $K_1$, $\tilde{K}_1$ and $K'_1$ are horizontal. Right: in the projection to the plane L.}
\label{fig2}
\end{center}
\end{figure}

\textit{Step 2.} By $s$ denote the boundary $(n-1)$-dimensional hypersphere of the half-sphere $C(x, \pi/2)$. Let $K$ be the intersection of the half-sphere $C(x, \pi/2)$ with the half-sphere bounded by $l$ and containing $x$. The spherical polytope $\psi_x \subset C(x, \pi/2)$ and $l$ contains its facet. Therefore, $\psi_x \subset K$. By $K'$ denote $f_{x, \lambda} (K)$ and by $\tilde{K}$ denote the locus of points inside $K$ that are at the distance not less than $\phi_1 - \alpha_1$ from $\partial K$ (see Figure~\ref{fig2}). Clearly, $\psi'_x \subset K'$ and it suffices to show that $K' \subset \tilde{K}$. 

The boundary of $K$ consists of two parts: $K_1 = K \cap l$ and $K_2 = K \cap s$. We divide in a similar way the boundaries of $K'$ and $\tilde{K}$: $\partial K' = K'_1 \cup K'_2$ and $\partial \tilde{K} = \tilde{K}_1 \cup \tilde{K}_2$.

For similarity, we use the same notations as in Step 1. Take an arbitrary point $a' \in \partial K'$. Let $\tilde{a} \in \tilde{K}$ be the intersection point of the arc-ray $xa'$ with $\tilde{K}$ and $a$ be the intersection point of this arc-ray with $l$. We should prove that $\tilde{a}$ lies between $a$ and $a'$. If $a' \in K'_1$ and $\tilde{a} \in \tilde{K}_1$ then in Step 1 we already proved the desired. If $a' \in K'_2$ and $\tilde{a} \in \tilde{K}_2$, then this follows by the same argument. If $a' \in K'_1$, then $\tilde{a}$ can not belong to $\tilde{K}_2$. Indeed, we can consider the $(n-1)$-dimensional hypersphere $h$ through $x$ and $l \cap s$. It divides $\partial K$ precisely into $K_1$ and $K_2$ and $\partial K'$ into $K'_1$ and $K'_2$. But $K$ is a part of sphere corresponding to an obtuse angle between $l$ and $s$. Therefore, the intersection of $\tilde{K}_1$ and $\tilde{K}_2$ is always from the same side with respect to $h$.

We reduced to prove our statement for $a' \in K'_2$ and $\tilde{a} \in \tilde{K}_1$. In this case $a$ lies outside $C(x, \pi/2)$ (and $a'$ does not coincide with $f_{x, \lambda}(a)$). Figure~\ref{fig1} is almost appropriate in this situation: we just should take into account that the angle $aOx$ is obtuse now.

We use the same notations for angles too: $\alpha_2 = a'Ox$, $\tilde{\alpha}_2 = \tilde{a}Ox$ and $\phi_2 = aOx$. Clearly, $\phi_2 > \pi/2$, but $\phi_2 - \tilde{\alpha}_2 \leq \pi/2$. The angle $\alpha$ does not depend on the choice of $a' \in K'_2$ and is equal to $\arcsin{\lambda}$. For the angle $\tilde{\alpha}_2$ equation (\ref{par}) still holds: $$\frac{\sin (\phi_2-\tilde{\alpha}_2)}{\sin\phi_2} = \frac{\sin(\phi_1-\alpha_1)}{\sin{\phi_1}} = c < 1.$$

Therefore, the angle $\tilde{\alpha}_2$ can be computed by the formula: $$\tilde{\alpha}_2 = \phi_2 - \arcsin(c\sin\phi_2).$$ In the beginning of the proof of Proposition 1 we already computed the derivative of this expression as a function of $\varphi$. We can see that it grows over the interval $(0; \pi)$. We can consider the case when $a'$ lies in the intersection of $K'_1$ and $K'_2$. For this case Step 1 implies that $\tilde{a}$ lies between $a$ and $a'$. But for any other case the angle $\tilde{\alpha}_2$ is bigger and the angle $\alpha_2$ is constant. This finishes the proof.

\section{Proof of Theorem 2}

We briefly sketch the main idea. If we perturb slightly the construction of the set $\Psi'$, then we can guarantee that any two points at the distance close to 1 have different colors. Then we can extend our coloring from the sphere to a small spherical shell and continue. The number of shells will not depend on $n$.

More precisely, fix some $R > 1/2$ and small $\e >0$. Consider a function $$g(r) = \frac{1}{2} - \frac{1}{4r^2} + \sqrt{\frac{1}{4}-\frac{5r^2-1}{16r^4}}.$$ Determine $r_{\ast} \in (\frac{\sqrt 5}{2}; + \infty)$ by the equation $g(r_{\ast}) = \frac{1}{2} + \e$. Let $\phi(r)$ be a continuous function
$$\phi(r) = \left\{
\begin{aligned}
&\arccos \sqrt{g(r)}, &\qquad r > r_{\ast}& \\
&\arccos \sqrt{g(r_{\ast})}, &\qquad \frac{1}{2} \leq r \leq r_{\ast}&
\end{aligned}
\right. $$

Let $\lambda(r) = \frac{1}{1+8\cos^2 \phi(r) + \e}$. Note that $\lambda(r)$ is a continuous function. It is obvious that for all $\e > 0$ we have $0 < \phi(r) < \pi/4$ and $0 < \lambda(r) < 1$.

Define $$\delta (r) = \min\{1-2r\lambda(r) \sin\big(2\phi(r)\big), 2r\sin\Big[\phi(r)-\arcsin\big(\lambda(r)\sin\phi(r)\big)\Big] - 1 \}.$$

Both expressions in the definition of $\delta$ are obtained from the system (\ref{se}). For any $r \in [1/2, +\infty)$ we are able to construct a set $\Psi'(r) \subset S^n_r$ with parameters $\phi(r)$ and $\lambda(r)$ (now we allow $r$ to be equal 1/2: our construction in the proof of Theorem 1 admits it). It can be easily obtained from Subsections 4.1 and 4.4 that the distance between any two points from this set is not in the interval $(1 - \delta(r); 1 + \delta(r))$. Therefore, an open $\delta/2$-neighborhood in $\R^{n+1}$ of $\Psi'(r)$ does not contain a pair of points at the distance 1.

Define $R_1 = R$ and $R_2 = R_1-\frac{\delta(R_1)}{2}$. We can obtain a proper coloring of the sphere $S^n_{R_1}$ using no more than $(\lambda^{-1}(R_1)+o(1))^n$ colors. Consider the spherical shell  $D_1 = \{x\in B^{n+1}_{R_1} | x \not\in B^{n+1}_{R_2} \}$. Extend our coloring to $D_1$. Color a point $x \in D_1$ with the color of the intersection point of the ray $Ox$ and the sphere $S^n_{R_1}$. We get the coloring of $D_1$ without distance 1. Now we make an independent coloring of the sphere $S^n_{R_2}$ (all colors are supposed to be different from the previous coloring) and extend it to the next spherical shell $D_2 = \{x\in B^{n+1}_{R_2} | x \not\in B^{n+1}_{R_3} \}$, where $R_3 = R_2-\frac{\delta(R_2)}{2}$. Then we continue in the same way.

The function $\delta(r)$ is continuous over the segment $[1/2; R]$. Therefore, it reaches its minimal value and this value is greater than 0. We can see that after a finite number of steps we get $R_k < 1/2$. If $R_k \leq 0$, we are done. If $0 < R_k < 1/2$, then it remains to color the ball $B^{n+1}_{R_k}$ in one color. It is clear that $\lambda(r)$ decreases monotonically as $r$ tends to 1/2. Therefore, we used no more than $k(\lambda^{-1}(R)+o(1))^n$ colors. Here $\lambda^{-1}(R)$ is equal to $x(R) + \e$, $k$ does not depend on $n$ and $\e$ is arbitrary. Hence, we get $$\chi(B^{n+1}_R) \leq (x(R)+o(1))^n.$$

\vskip+0.2cm
\textbf{Acknowledgements.} The author is grateful to A.M. Raigorodskii, N.G. Moshchevitin, A.B. Kupavskii, D. Cherkashin and A. Polyanskii for their interest in this work and fruitful discussions and the anonymous referees for useful remarks.
\vskip+0.2cm

\bibliographystyle{abbrv}
\bibliography{sphere_chromatic}

\begin{thebibliography}{10}

\bibitem{AR}
S.~Artstein-Avidan and O.~Raz.
\newblock Weighted covering numbers of convex sets.
\newblock {\em Adv. Math.}, 227(1):730--744, 2011.

\bibitem{AS}
S.~Artstein-Avidan and B.~A. Slomka.
\newblock On weighted covering numbers and the {L}evi-{H}adwiger conjecture.
\newblock {\em Israel J. Math.}, 209(1):125--155, 2015.

\bibitem{BPT}
C.~Bachoc, A.~Passuello, and A.~Thiery.
\newblock The density of sets avoiding distance 1 in {E}uclidean space.
\newblock {\em Discrete Comput. Geom.}, 53(4):783--808, 2015.

\bibitem{Be}
A.~V. Berdnikov.
\newblock Estimate for the chromatic number of {E}uclidean space with several
  forbidden distances.
\newblock {\em Mat. Zametki}, 99(5):783--787, 2016.

\bibitem{Be2}
A.~V. Berdnikov.
\newblock Chromatic numbers of distance graphs with several forbidden distances
  and without cliques of a given size.
\newblock {\em Problems of Information Transmission}, 54(1):70--83, 2018.

\bibitem{BeRa}
A.~V. Berdnikov and A.~M. Raigorodskii.
\newblock On the chromatic number of {E}uclidean space with two forbidden
  distances.
\newblock {\em Math. Notes}, 96(5-6):827--830, 2014.
\newblock Translation of Mat. Zametki {{\bf{9}}6} (2014), no. 5, 790--793.

\bibitem{Bo}
K.~B\"or\"oczky, Jr.
\newblock {\em Finite packing and covering}, volume 154 of {\em Cambridge
  Tracts in Mathematics}.
\newblock Cambridge University Press, Cambridge, 2004.

\bibitem{BW}
K.~B\"or\"oczky, Jr. and G.~Wintsche.
\newblock Covering the sphere by equal spherical balls.
\newblock In {\em Discrete and computational geometry}, volume~25 of {\em
  Algorithms Combin.}, pages 235--251. Springer, Berlin, 2003.

\bibitem{CR}
D.~D. Cherkashin and A.~M. Raigorodskii.
\newblock On the chromatic numbers of small-dimensional spaces.
\newblock {\em Dokl. Math.}, 95(1):5--6, 2017.

\bibitem{EdB}
N.~G. de~Bruijn and P.~Erd\"os.
\newblock A colour problem for infinite graphs and a problem in the theory of
  relations.
\newblock {\em Nederl. Akad. Wetensch. Proc. Ser. A. {\bf 54} = Indagationes
  Math.}, 13:369--373, 1951.

\bibitem{dG}
A.~D.~N.~J. {de Grey}.
\newblock {The chromatic number of the plane is at least 5}.
\newblock {\em ArXiv e-prints}, Apr. 2018.

\bibitem{DfW}
E.~{DeCorte}, F.~{M{\'a}rio de Oliveira Filho}, and F.~{Vallentin}.
\newblock {Complete positivity and distance-avoiding sets}.
\newblock {\em ArXiv e-prints}, Apr. 2018.

\bibitem{Fa}
K.~J. Falconer.
\newblock The realization of distances in measurable subsets covering {${\bf
  R}^{n}$}.
\newblock {\em J. Combin. Theory Ser. A}, 31(2):184--189, 1981.

\bibitem{Jo}
D.~S. Johnson.
\newblock Approximation algorithms for combinatorial problems.
\newblock {\em J. Comput. System Sci.}, 9:256--278, 1974.
\newblock Fifth Annual ACM Symposium on the Theory of Computing (Austin, Tex.,
  1973).

\bibitem{Ko}
O.~A. Kostina.
\newblock On lower bounds for the chromatic number of a sphere.
\newblock {\em to appear in Math. Notes}.

\bibitem{KRa}
O.~A. Kostina and A.~M. Raigorodskii.
\newblock On lower bounds for the chromatic number of a sphere.
\newblock {\em Dokl. Math.}, 92(1):500--502, 2015.

\bibitem{Ku1}
A.~Kupavskiy.
\newblock The chromatic number of the space $\mathbb{R}^n$ with the set of
  forbidden distances.
\newblock {\em Dokl. Math.}, 82(3):963--966, 2010.

\bibitem{Ku2}
A.~Kupavskiy.
\newblock On the chromatic number of {$\Bbb R^n$} with an arbitrary norm.
\newblock {\em Discrete Math.}, 311(6):437--440, 2011.

\bibitem{Lar}
D.~G. Larman and C.~A. Rogers.
\newblock The realization of distances within sets in {E}uclidean space.
\newblock {\em Mathematika}, 19:1--24, 1972.

\bibitem{Lo}
L.~Lov\'asz.
\newblock On the ratio of optimal integral and fractional covers.
\newblock {\em Discrete Math.}, 13(4):383--390, 1975.

\bibitem{Lo1}
L.~Lov\'asz.
\newblock Self-dual polytopes and the chromatic number of distance graphs on
  the sphere.
\newblock {\em Acta Sci. Math. (Szeged)}, 45(1-4):317--323, 1983.

\bibitem{Na}
M.~Nasz\'odi.
\newblock On some covering problems in geometry.
\newblock {\em Proc. Amer. Math. Soc.}, 144(8):3555--3562, 2016.

\bibitem{NaS}
M.~{Nasz{\'o}di}.
\newblock Flavors of translative coverings.
\newblock In {\em New Trends in Intuitive Geometry}. Springer, Berlin, 2018.

\bibitem{NP}
M.~Nasz\'odi and A.~Polyanskii.
\newblock Approximating set multi-covers.
\newblock {\em European J. Combin.}, 67:174--180, 2018.

\bibitem{PoRa}
E.~I. Ponomarenko and A.~M. Raigorodskii.
\newblock New lower bound for the chromatic number of a rational space with one
  and two forbidden distances.
\newblock {\em Mat. Zametki}, 97(2):255--261, 2015.

\bibitem{Pr1}
R.~{Prosanov}.
\newblock {A new proof of the Larman-Rogers upper bound for the chromatic
  number of the Euclidean space}.
\newblock {\em ArXiv e-prints}, Oct. 2016.

\bibitem{Pr2}
R.~I. Prosanov.
\newblock Upper bounds for the chromatic numbers of {E}uclidean spaces with
  forbidden {R}amsey sets.
\newblock {\em Mat. Zametki}, 103(2):248--257, 2018.

\bibitem{Ra1}
A.~M. Raigorodskii.
\newblock On the chromatic number of a space.
\newblock {\em Uspekhi Mat. Nauk}, 55(2(332)):147--148, 2000.

\bibitem{Ra2}
A.~M. Raigorodskii.
\newblock The {B}orsuk problem and the chromatic numbers of some metric spaces.
\newblock {\em Uspekhi Mat. Nauk}, 56(1(337)):107--146, 2001.

\bibitem{Rai1}
A.~M. Raigorodskii.
\newblock On the chromatic numbers of spheres in {E}uclidean spaces.
\newblock {\em Dokl. Akad. Nauk}, 432(2):174--177, 2010.

\bibitem{Rai}
A.~M. Raigorodskii.
\newblock On the chromatic numbers of spheres in {${\Bbb R}^n$}.
\newblock {\em Combinatorica}, 32(1):111--123, 2012.

\bibitem{Rai2}
A.~M. Raigorodskii.
\newblock Coloring distance graphs and graphs of diameters.
\newblock In {\em Thirty essays on geometric graph theory}, pages 429--460.
  Springer, New York, 2013.

\bibitem{Rai3}
A.~M. Raigorodskii.
\newblock Cliques and cycles in distance graphs and graphs of diameters.
\newblock In {\em Discrete geometry and algebraic combinatorics}, volume 625 of
  {\em Contemp. Math.}, pages 93--109. Amer. Math. Soc., Providence, RI, 2014.

\bibitem{Ra}
A.~M. Raigorodskii.
\newblock Combinatorial geometry and coding theory.
\newblock {\em Fund. Inform.}, 145(3):359--369, 2016.

\bibitem{Ro1}
C.~A. Rogers.
\newblock A note on coverings.
\newblock {\em Mathematika}, 4:1--6, 1957.

\bibitem{Ro2}
C.~A. Rogers.
\newblock Covering a sphere with spheres.
\newblock {\em Mathematika}, 10:157--164, 1963.

\bibitem{Sa2}
A.~Sagdeev.
\newblock On a {F}rankl--{R}\"{o}dl theorem and its geometric corollaries.
\newblock {\em Electronic Notes in Discrete Mathematics}, 61:1033 -- 1037,
  2017.

\bibitem{Sa}
A.~A. Sagdeev and A.~M. Raigorodskii.
\newblock On the chromatic number of a space with forbidden regular simplices.
\newblock {\em Dokl. Math.}, 95(1):15--16, 2017.

\bibitem{So}
A.~Soifer.
\newblock {\em The mathematical coloring book}.
\newblock Springer, New York, 2009.
\newblock Mathematics of coloring and the colorful life of its creators.

\bibitem{St}
S.~K. Stein.
\newblock Two combinatorial covering theorems.
\newblock {\em J. Combinatorial Theory Ser. A}, 16:391--397, 1974.

\end{thebibliography}

\end{document}